\documentclass[a4paper,12pt]{article}
\usepackage{amsmath,amssymb}
\usepackage{enumerate}
\usepackage{epic}
\usepackage{bbm}
\usepackage{theorem}

\usepackage{bbm}
\usepackage{mathptmx}

\newcommand{\NN}{\mathbb{N}}
\newcommand{\ZZ}{\mathbb{Z}}

\newcommand{\RR}{\mathbb{R}}

\newcommand{\cF}{\ensuremath{\mathcal{F}}}
\newcommand{\cL}{\mathcal{L}}
\newcommand{\cB}{\mathcal{B}}

\newcommand{\<}{\ensuremath{\langle}}
\renewcommand{\>}{\ensuremath{\rangle}}

\newcommand{\supp}{\operatorname{supp}}
\newcommand{\ini}{\operatorname{in}}
\newcommand{\bin}{\operatorname{bin}}
\newcommand{\mon}{\operatorname{mon}}
\newcommand{\sat}{\operatorname{sat}}

\newcommand{\F}{\ensuremath{\mathcal{F}}}

\newcommand{\face}{\operatorname{face}}
\newcommand{\conv}{\operatorname{conv}}
\newcommand{\Span}{\operatorname{span}}
\newcommand{\State}{\operatorname{State}}

\renewcommand{\>}{\rangle}

\newtheorem{Theorem}{Theorem}[section]
\newtheorem{Algorithm}[Theorem]{Algorithm}
\newtheorem{Lemma}[Theorem]{Lemma}
\newtheorem{Proposition}[Theorem]{Proposition}
\newtheorem{Corollary}[Theorem]{Corollary}

\newtheorem{Example}[Theorem]{Example}

\newtheorem{Definition}[Theorem]{Definition}

\newenvironment{Proof}{\noindent{\it Proof.\/}}{\hfill $\square$\medskip}

\begin{document}
\title{Truncated Gr\"obner fans and lattice ideals}
\author{
Niels Lauritzen\\
Institut for Matematiske Fag\\
Aarhus Universitet\\
DK-8000 \AA rhus\\
Denmark\\ 
niels@imf.au.dk
}

\maketitle

\begin{abstract}
We outline a generalization of the Gr\"obner fan of
a homogeneous ideal with maximal cells parametrizing 
truncated Gr\"obner bases. This ``truncated'' 
Gr\"obner fan is usually much smaller 
than the full Gr\"obner fan and offers the natural framework
for conversion between truncated Gr\"obner bases. The 
generic Gr\"obner walk generalizes naturally 
to this setting by using the Buchberger algorithm
with truncation on facets. 

We specialize to the setting of lattice ideals.
Here facets along the generic walk are given by
unique (facet) binomials. This along with the 
representation of binomials as integer vectors
give an especially simple version of the
generic Gr\"obner walk.

Computational experience with the special
Aardal-Lenstra integer programming knapsack problems is
reported.
The algorithms of this paper are implemented in the 
software package {\tt GLATWALK}, which is available for download at 
{\tt http://home.imf.au.dk/niels/GLATWALK}.
\end{abstract}


\section{Introduction}

The
generic Gr\"obner walk \cite{FJLT} is a version of the classical Gr\"obner
walk algorithm for Gr\"obner basis conversion in the Gr\"obner
fan of an ideal in a polynomial ring. 
In the generic walk explicit rational vectors in
the Gr\"obner fan are replaced by computations with infinitesimal
numbers which can be handled formally. This leads to an algorithm
where input consists only of a source Gr\"obner basis, a source
term order and a target term order.

Truncation of homogeneous ideals have proved very valuable in algebraic
computations related to integer programming (see for example
\cite{TW}). In general there are a lot fewer truncated initial
ideals than initial ideals. Similary to initial ideals, 
truncated initial ideals may be parametrized by the maximal
cells in a complete polyhedral fan. We introduce this
fan, which is easily constructed
from the usual Gr\"obner fan by inserting the truncation operator
at the appropriate places.
The {\it truncated Gr\"obner fan\/}
is in general much smaller than the full Gr\"obner fan and forms
the polyhedral setting for a truncated version of the
generic Gr\"obner walk. We prove that 
the truncated Gr\"obner fan is regular 
along the lines of \cite{St}. This leads to a ``truncated''
state polytope with vertices enumerating the different reduced
truncated Gr\"obner bases.

In the setting of lattice ideals we give
a rather detailed version of the generic Gr\"obner walk.
Algebraic
computations with lattice ideals can be greatly simplified
representing (saturated) binomials as integer vectors. This along
with the fact that the generic walk only traverses facets
lead to several simplifications.
We report on computational experience in computing saturations
of lattice ideals and truncated test sets related to the
integer programming problems posed in \cite{LATTE}. Our experiments
show that the the generic walk in the truncated Gr\"obner fan 
consists of significantly fewer steps, whereas the walk in the
full Gr\"obner fan does not compare well computing directly
with the Buchberger algorithm.

I am grateful to B.~Sturmfels for inspiring conversation and
for greatly simplifying my original approach to truncated
Gr\"obner fans. Thanks are also due to 
K.~Fukuda, A.~N.~Jensen and R.~Thomas for slowly making me
grasp the joys of algebra and polyhedral geometry.

\section{Preliminaries}

We let $R = k[x_1, \dots, x_n]$ denote the ring of polynomials
over a field $k$. We will view $R$ as the semigroup ring
$k[\NN^n]$. 

\subsection{Grading on $R$}
Given $n$ elements $a_1, \dots, a_n$ of
an abelian group $(A, +)$ we let
$S_A$ denote the semigroup $\NN a_1 + \cdots + \NN a_n\subset A$. For 
$v = (v_1, \dots, v_n)\in \NN^n$ we put 
$\deg(x^v) = v\cdot a = v_1 a_1 + \cdots + v_n a_n$. These data give a natural
$A$-grading on $R$ by defining
$$
R_s = \Span_k \{x^v \mid \deg(x^v) = s\}
$$  
for $s\in S_A$. 
A non-zero element $f\in R$ is called
homogeneous of degree $\deg(f) = s$ if $f\in R_s$.
Given an ideal $J$ in $R$, we let $J_s = J\cap R_s$.
Recall that an ideal $J$ is homogeneous ideal if 
$J = \oplus J_s$ and that this is
equivalent to $J$ being generated by homogeneous elements.
We call $R$ {\it positively graded\/} if $R_0 = k$. This
is equivalent to $\dim_k R_s < \infty$ for every $s\in S_A$.

\subsection{Truncating subsets}
A subset $\Omega\subset S_A$ 
is called {\it truncating\/} if $s, t\in \Omega$ whenever $s+t\in \Omega$ for
$s, t\in S_A$. 
Our standard example of a truncating subset is
$$
\Omega_b = \{x\in S_A \mid b - x\in S_A\}
$$
for $b\in S_A$. To a truncating subset $\Omega \subset S_A$ we
associate the monomial ideal (cf.~Remark 12.8 in \cite{St})
$$
M_\Omega = \<x^v \mid \deg(v)\not\in \Omega\>\subset R.
$$
Given a homogeneous ideal $J \subset R$ we let
$$
J_\Omega = \bigoplus_{s\in \Omega} J_s.
$$
\begin{Lemma}\label{LemmaTE}
Let $I$ and $J$ be homogeneous ideals in $R$. Then
$$
I_\Omega = J_\Omega \mathrm{\ if\ and\ only\ if\ } 
I + M_\Omega = J + M_\Omega.
$$
\end{Lemma}
\begin{Proof}
If $I$ is a homogeneous ideal, then $I + M_\Omega$ is 
a homogeneous ideal and
$$
(I + M_\Omega)_s = 
\left\{
\begin{array}{ll}
I_s & \mbox{if $s\in \Omega$}\\
(M_\Omega)_s & \mbox{if $s\not\in\Omega$}.
\end{array}
\right.
$$
This proves the lemma.

\end{Proof}

\subsection{Initial ideals}

Let $\prec$ denote 
a total multiplicative ordering on monomials in $R$ (we do 
not require that the monomial $1$ is minimal). If $\omega\in \RR^n$, we let 
$\prec_\omega$ denote the multiplicative ordering defined by 
$x^u \prec_\omega x^v$ if $\omega\cdot u < \omega \cdot v$ or 
$\omega\cdot u = \omega\cdot v$ and $x^u \prec x^v$.
For $f\in R$, we let $\supp(f)$ denote the set of monomials occuring
with non-zero coefficient in $f$. We let $\ini_\prec(f)$ denote
the maximal (initial) term in $\supp(f)$ with respect to $\prec$. Similarly we
let $\ini_\omega(f)$ denote the sum of terms $a_v x^v$ in $\supp(f)$
with $\omega\cdot v$ maximal.
For a subset $G\subset R$
we let  $\ini_\omega(G)$ and  $\ini_\prec(G)$ denote the ideals 
$\<\ini_\omega(f)\mid f\in G\setminus\{0\}\>$ 
and $\<\ini_\prec(f) \mid f\in G\>$ respectively. These
ideals are homogeneous if $I$ is homogeneous. 
A {\it Gr\"obner basis} for $I$ over $\prec$ is a finite set
$G:=\{f_1, \dots, f_r\} \subset I$ such that
$$
\ini_\prec(I) = \ini_\prec(G) = \<\ini_\prec(f_1), \dots, \ini_\prec(f_r)\>.
$$
The Gr\"obner basis $G$ is called {\it minimal\/} if none of $f_1, \dots, 
f_r$ can be left out and {\it reduced\/} if $\ini_\prec(f_i)$ does
not divide any of the terms in $f_j$ for $i\neq j$ and $i, j = 1, \dots, r$.
The reduced Gr\"obner basis of an ideal is unique and consists of
homogeneous elements if the ideal is homogeneous.
A homogeneous ideal $J$ in $R$
always has a reduced Gr\"obner basis over $\prec$ if
$R$ is positively graded, since $\dim_k J_s < \infty$ for
$s\in S_A$.
We record the following simple but crucial result (\cite{St}, Proposition 1.8) 
with a complete proof.

\begin{Proposition}
\label{PropositionCrux}
Let $I\subset R$ be any ideal and $\omega\in \RR^n$. Then
$$
\ini_\prec(\ini_\omega(I)) = \ini_{\prec_\omega}(I).
$$
\end{Proposition} 
\begin{Proof}
Clearly $\ini_{\prec_\omega}(I) \subset \ini_\prec(\ini_\omega(I))$. The
ideal $\ini_\omega(I)$ is homogeneous in the grading given by $\omega$. So
we may decompose an element $f\in \ini_\omega(I)$ as 
$f = f_{\lambda_1} + \cdots + f_{\lambda_t}$, where $f_{\lambda_j}\in
\ini_\omega(I)$ is homogeneous of $\omega$-weight $\lambda_j$ for
$j = 1, \dots, t$.
Now $\ini_\prec(f) = 
\ini_\prec(f_{\lambda_j})$ for some $j$ and we may write
$$
f_{\lambda_j} = a_1\ini_\omega(f_1) + \cdots + a_r \ini_\omega(f_r)
$$
for suitable $f_1, \dots, f_r\in I$, where $a_1, \dots, a_r$ are
homogeneous elements. Therefore 
$$
f_{\lambda_j} = \ini_\omega(a_1 f_1 + \cdots + a_r f_r)
$$
and $\ini_\prec(f) = \ini_{\prec_\omega}(a_1 f_1 + \cdots + a_r f_r)$.
This shows that $\ini_{\prec_\omega}(I)\supset \ini_\prec(\ini_\omega(I))$.
\end{Proof}

Notice that the multiplicativity of $\prec$ is not used in the proof
of Proposition \ref{PropositionCrux}. The lifting from $\ini_\omega(I)$
to $I$ in the proof is a key element in the Gr\"obner
walk algorithm.

\subsection{Truncated Gr\"obner bases}

Let $J$ be a homogeneous ideal in $R$ and $\Omega$ a truncating
subset.
A finite subset $G\subset J$ is called an
$\Omega$-Gr\"obner basis for $J$ over $\prec$ if
$$
\ini_\prec(J)_\Omega = \ini_\prec(G)_\Omega.
$$
If the coefficients of the initial terms $\ini_\prec(g)$ are
$1$ for $g\in G$ and $\ini_\prec(g)$ does not divide any of
the terms in $g'$ for $g\neq g'\in G$, then $G$ is called
a reduced $\Omega$-Gr\"obner basis. Reduced $\Omega$-Gr\"obner
bases are unique.

\begin{Proposition}
\label{PropositionRG}
Let $G$ be the reduced Gr\"obner basis for $J$ over $\prec$.
Then $G$ consists of homogeneous elements and 
$G_\Omega = \{g\in G\mid \deg(g) \in \Omega\}$ is the reduced
$\Omega$-Gr\"obner basis for $J$ over $\prec$.
\end{Proposition}
\begin{Proof}
If $g$ is an element of the reduced Gr\"obner basis of $J$ and
$g = g_1 + \cdots + g_r$ is written as a sum of homogeneous elements,
then $\ini_\prec(g) = \ini_\prec(g_j)$ for some $j = 1, \dots, r$. Therefore
$g$ has to be homogeneous.
The monomial ideal $\ini_\prec(J)$ is spanned as a 
vector space by $\{\ini_\prec(f) \mid f\in J_s, \mathrm{for\ some\ }s\in S\}$.
This shows that $\ini_\prec(J)_\Omega$ is the $k$-span of $\ini_\prec(f)$ for
$f\in J_\Omega$. Suppose that $f\in J_\Omega$. 
Since $G$ is a Gr\"obner basis for $J$ we may find
$g\in G$, such that $\ini_\prec(g)$ divides $\ini_\prec(f)$. This shows that
$\deg(g) \in \Omega$, since $\Omega$ is a truncating subset. Therefore
$\ini_\prec(f)\in \<\ini_\prec(g)\mid g\in G_\Omega\>_\Omega$.
\end{Proof}

\begin{Corollary}\label{CorollaryDeformed}
If $G_\Omega$ is the reduced 
$\Omega$-Gr\"obner basis for $I$ over $\prec_\omega$ 
then  $\{\ini_\omega(g) \mid g\in G_\Omega\}$ is the reduced 
$\Omega$-Gr\"obner basis for $\ini_\omega(I)$ over $\prec$
\end{Corollary}
\begin{Proof}
Let $G$ be the reduced Gr\"obner basis for $I$ over $\prec_\omega$.
Then we know from Proposition \ref{PropositionCrux} that 
$\{\ini_\omega(g) \mid g\in G\}$ is the reduced Gr\"obner basis of
$\ini_\omega(I)$. Now Proposition \ref{PropositionRG} gives that
$$
\{\ini_\omega(g) \mid g\in G, \deg(\ini_\omega(g))\in \Omega\} 
=\{\ini_\omega(g) \mid g\in G_\Omega\}
$$
is the reduced $\Omega$-Gr\"obner basis of $I$.
\end{Proof}

Proposition \ref{PropositionRG} reveals that the
$\Omega$-truncated Gr\"obner basis can be obtained
from the reduced Gr\"obner basis by picking out the
elements with degree in $\Omega$. Truncated Gr\"obner bases
can be computed from a homogeneous generating set using Buchbergers 
algorithm discarding $S$-polynomials with degree
outside $\Omega$. This follows from the fact that the
division algorithm preserves the degree of a homogeneous polynomial.
Very often only Gr\"obner bases up to a certain degree are needed.

\section{The truncated Gr\"obner fan}\label{SectionFan}

In this section we assume that $R$ is positively graded. 
Analogously to (\cite{St}, Proposition 2.3) we define 
$$
C_\Omega[\omega] = \{v\in \RR^n \mid \ini_v(I)_\Omega =
\ini_\omega(I)_\Omega\}
$$
for a homogeneous ideal $I\subset R$. We call the 
closure of $C_\Omega[\omega]$ in $\RR^n$ a {\it truncated Gr\"obner cone\/}.

\begin{Theorem}
The collection 
$$
\F_\Omega(I) = \{\overline{C_\Omega[\omega]} \mid\omega\in \RR^n\}
$$
of truncated Gr\"obner cones form a complete fan in $\RR^n$.
\end{Theorem}
\begin{Proof}
A monomial ideal $M$ satisfies $\ini_u(J) + M = 
\ini_u(J + M)$ for every ideal $J\subset R$ and $u\in
\RR^n$. Now
Lemma \ref{LemmaTE} shows that $\ini_v(I)_\Omega = \ini_\omega(I)_\Omega$ if
and only if $\ini_v(I + M_\Omega) = \ini_\omega(I + M_\Omega)$. This proves
that $\F_\Omega(I)$ is the usual Gr\"obner fan of the homogeneous 
ideal $I + M_\Omega$. Now the conclusion follows from
(\cite{St}, Proposition 2.4).
\end{Proof}

The truncated Gr\"obner fan is available from the usual Gr\"obner
fan by eliminating inequalities given by polynomials of
degree outside $\Omega$. We give an example illustrating this.

\begin{Example}

Consider the (toric) ideal 
$$
I_A = \<a^2 c - b^2 e, a^2 d - b e^2, c e - b d\>\subset k[a, b, c, d, e].
$$
This ideal is homogeneous in the grading given by the columns of  
$$
A =
\begin{pmatrix}
1 & 1 & 1 & 1 & 1\\
0 & 1 & 2 & 1 & 0\\
0 & 0 & 1 & 2 & 1
\end{pmatrix}.
$$
For example, the degree of the variable $c$ is $(1, 2, 1)$.
The
Gr\"obner fan $\F(I_A)$ of $I_A$ is the normal fan of
an octagon in $\RR^5$ (cf.~Example 1.1 in \cite{HT}). It is pictured 
in $\RR^2$ below with reduced Gr\"obner bases labeling the maximal cells.
$$
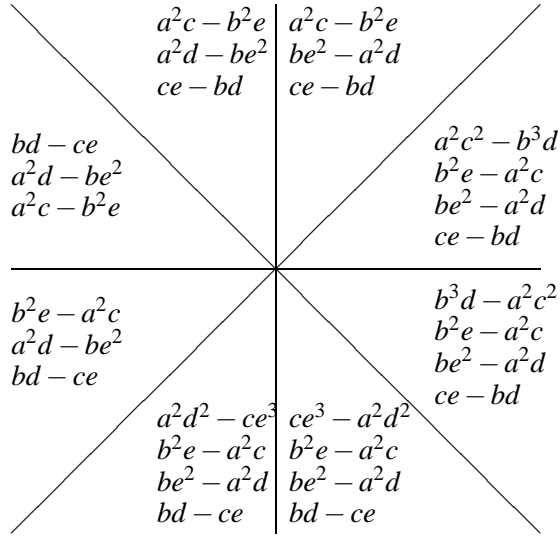
\begin{picture}(200,200)(0,0)

\put(0,100){\line(1,0){200}}
\put(0,0){\line(1,1){200}}
\put(0,200){\line(1,-1){200}}
\put(100,0){\line(0,1){200}}


\put(55,40){{\small $a^2 d^2 - c e^3$}}
\put(55,28){{\small $b^2 e - a^2 c$}}
\put(55,16){{\small $b e^2 - a^2 d$}}
\put(55,4){{\small $b d - c e$}}


\put(105,40){{\small $c e^3 - a^2 d^2$}}
\put(105,28){{\small $b^2 e - a^2 c$}}
\put(105,16){{\small $b e^2 - a^2 d$}}
\put(105,4){{\small $b d - c e$}}


\put(0,80){{\small $b^2 e - a^2 c$}}
\put(0,68){{\small $a^2 d - b e^2$}}
\put(0,56){{\small $b d - c e$}}


\put(0,120){{\small $a^2 c - b^2 e$}}
\put(0,132){{\small $a^2 d - b e^2$}}
\put(0,144){{\small $b d - c e$}}


\put(55,190){{\small $a^2 c - b^2 e$}}
\put(55,178){{\small $a^2 d - b e^2$}}
\put(55,166){{\small $c e - b d$}}


\put(105,190){{\small $a^2 c - b^2 e$}}
\put(105,178){{\small $b e^2 - a^2 d$}}
\put(105,166){{\small $c e - b d$}}


\put(160,145){{\small $a^2 c^2 - b^3 d$}}
\put(160,133){{\small $b^2 e - a^2 c$}}
\put(160,121){{\small $b e^2 - a^2 d$}}
\put(160,109){{\small $c e - b d$}}


\put(160,85){{\small $b^3 d - a^2 c^2$}}
\put(160,73){{\small $b^2 e - a^2 c$}}
\put(160,61){{\small $b e^2 - a^2 d$}}
\put(160,49){{\small $c e - b d$}}

\end{picture}
$$
Putting $\Omega = \{v\in S_{\ZZ^3}\mid v\cdot (1,1,1) < 6\}$ we get
the truncated fan $\F_\Omega(I_A)$ with the reduced $\Omega$-Gr\"obner bases
labeling the maximal cells.
$$
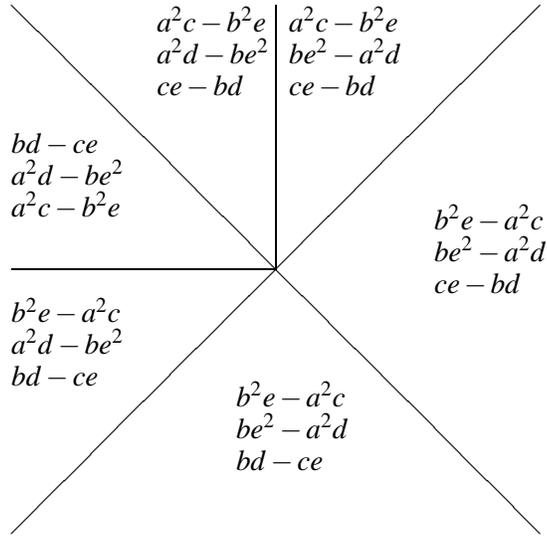
\begin{picture}(200,200)(0,0)

\put(0,100){\line(1,0){100}}
\put(0,0){\line(1,1){200}}
\put(0,200){\line(1,-1){200}}
\put(100,100){\line(0,1){100}}


\put(85,48){{\small $b^2 e - a^2 c$}}
\put(85,36){{\small $b e^2 - a^2 d$}}
\put(85,24){{\small $b d - c e$}}


\put(0,80){{\small $b^2 e - a^2 c$}}
\put(0,68){{\small $a^2 d - b e^2$}}
\put(0,56){{\small $b d - c e$}}


\put(0,120){{\small $a^2 c - b^2 e$}}
\put(0,132){{\small $a^2 d - b e^2$}}
\put(0,144){{\small $b d - c e$}}


\put(55,190){{\small $a^2 c - b^2 e$}}
\put(55,178){{\small $a^2 d - b e^2$}}
\put(55,166){{\small $c e - b d$}}


\put(105,190){{\small $a^2 c - b^2 e$}}
\put(105,178){{\small $b e^2 - a^2 d$}}
\put(105,166){{\small $c e - b d$}}


\put(160,115){{\small $b^2 e - a^2 c$}}
\put(160,103){{\small $b e^2 - a^2 d$}}
\put(160,91){{\small $c e - b d$}}

\end{picture}
$$

\end{Example}

\section{Truncated state polytopes}\label{SectionState}
The truncated Gr\"obner fan of $I$ is the usual Gr\"obner fan
of $I + M_\Omega$. From this it follows that 
the truncated Gr\"obner fan
of a homogeneous ideal $I$ is 
the normal fan of a natural 
Minkowski summand in a state polytope for $I$. 
Emphasizing the simple Lemma \ref{LemmaRefine} below, we briefly sketch 
a proof of this along the lines of \cite{St}. 
Let $I$ be a homogeneous ideal in $R$ and $a\in A$.
\begin{Lemma}\label{LemmaRefine}
Let $\prec$ and $\prec'$ be two total multiplicative orderings on monomials
in $R$. Let
$x^{v_1}, \dots, x^{v_s}$ denote the monomials in $\ini_\prec(I)_a$
and $x^{u_1}, \dots, x^{u_s}$ the monomials in $\ini_{\prec'}(I)_a$.
Then after permuting $u_1, \dots, u_s$ and $v_1, \dots, v_s$ we may 
assume that $x^{v_1}\prec x^{v_2} \prec \cdots \prec x^{v_s}$ and
\begin{align*}
x^{u_1} &\prec x^{v_1}\\
x^{u_2} &\prec x^{v_2}\\
&\vdots\\ 
x^{u_s} &\prec x^{v_s}.
\end{align*}
\end{Lemma}
\begin{Proof}
We may find a vector space basis $f_1, \dots, f_s$ of $I_a$, such that
$\ini_\prec(f_1) = x^{v_1}, \dots, \ini_\prec(f_s) = x^{v_s}$ and 
$x^{v_1} \prec x^{v_2} \prec \cdots \prec x^{v_s}$.
Now put $f_1' = f_1$. Move on to $f_2$. If
$\ini_{\prec'}(f_2) = \lambda \ini_{\prec'}(f_1')$ for $\lambda\in k$, 
put $f_2' = f_2-\lambda f_1'$. Then clearly
$\ini_{\prec'}(f_2')\prec \ini_{\prec}(f_2)$. In general if
$\ini_{\prec'}(f_j')\prec \ini_{\prec}(f_j)$ for
$j < m$ and $\ini_{\prec'}(f_m)\in W = \Span_k \{\ini_{\prec'}(f_1'), \dots,
\ini_{\prec'}(f_{m-1}')\}$, then $f_m' = f_m - \lambda_1 f_1' - \cdots -
\lambda_{m-1} f_{m-1}'$ satisfies $\ini_{\prec'}(f_m')\not\in W$ for
suitable $\lambda_1, \dots, \lambda_{m-1}\in k$. Furthermore
$\ini_{\prec'}(f_m')\prec \ini_{\prec}(f_m)$. In this way we get
the monomials $x^{u_1} = \ini_{\prec'}(f_1'), \dots, 
x^{u_s} = \ini_\prec(f_s')$ of $\ini_{\prec'}(I)_a$ written up in
the desired way.
\end{Proof}

Let
$$
\Sigma\, \ini_\prec(I)_a = \sum_{x^v\in \ini_\prec(I)_a} v
$$
where $a\in A$ and $\prec$ is a multiplicative total ordering.

\begin{Corollary}\label{Cor1}
Let $\prec_1$ and $\prec_2$ be two total multiplicative orderings. If
$\Sigma\,\ini_{\prec_1}(I)_a = \Sigma\,\ini_{\prec_2}(I)_a$, then
$\ini_{\prec_1}(I)_a = \ini_{\prec_2}(I)_a$.
\end{Corollary}
\begin{Proof}
This is an easy consequence of Lemma \ref{LemmaRefine}.
\end{Proof}

\begin{Definition}
A state polytope in degree $a$ for $I$ is defined as
$$
\State_a(I) = \conv\{\Sigma\,\ini_\prec(I)_a \mid \prec\text{\ total multiplicative ordering}\}.
$$
For a finite subset $S\subset S_A$ we let
$$
\State_S(I) = \sum_{a\in S} \State_a(I).
$$
\end{Definition}

\begin{Corollary}\label{Cor2}
For $\omega\in \RR^n$ we have
$$
\face_\omega(\State_a(I)) = \State_a(\ini_\omega(I)).
$$
\end{Corollary}
\begin{Proof}
First assume that $\ini_\omega(I)$ is a
monomial ideal and $\face_\omega(\State_a(I))$ is a vertex
$\{\Sigma\,\ini_\prec(I)_a\}$ for some multiplicative monomial ordering $\prec$ 
by  picking a generic $\omega$. Apply Lemma \ref{LemmaRefine} to
$\ini_{\prec_\omega}(I)_a$ and $\ini_\prec(I)_a$. In this setting we then have
$u_1 = v_1, \dots, u_s = v_s$ since
$\omega \cdot (v_1 + \cdots + v_s) \geq \omega\cdot(u_1 + \cdots
+ u_s)$. Therefore $\{\Sigma\,\ini_\prec(I)_a\} = \{\Sigma\,\ini_\omega(I)_a\}$
when $\omega$ is generic. 
General $\omega$ are reduced to generic $\omega$ as in the last part
of the proof of Lemma 2.6 in \cite{St}. 
\end{Proof}

Now let $G$ be a universal Gr\"obner basis  for $I$ consisting of homogeneous
elements. Put $S = \{\deg(g)\mid g\in G\}$. 
Given a truncating subset $\Omega\subset S_A$ we let 
$$
\State_\Omega(I) = \State_{\Omega\cap S}(I)
$$
denote a {\it truncated state polytope\/}.
Notice that $\State_\Omega(I)$ is a Minkowski summand in a
state polytope $\State_S(I)$ of $I$. Using Corollaries \ref{Cor1} and
\ref{Cor2}, the same arguments as in 
the last part of the proof of Theorem 2.5 in \cite{St} show 
that the normal fan of $\State_\Omega(I)$ is
$\F_\Omega(I)$. It follows that $\F_\Omega(I)$ is the normal fan of a Minkowksi
summand of a state polytope and that $\F_\Omega(I)$ is
a coarsening of the usual Gr\"obner fan.

\section{Walking in the truncated Gr\"obner fan}\label{SectionWalk}

The Gr\"obner walk can be carried out in the truncated Gr\"obner
fan converting one truncated Gr\"obner basis to another. We sketch
the appropriate generalization of Proposition 3.2 in \cite{FJLT}. The
term orders of (\cite{FJLT}, Proposition 3.2) are represented by
weight vectors below. If a weight vector $\omega$ represents
the term order $\prec$, then we refer to $\prec_\eta$ as
$\omega$ modified by $\eta$.

\begin{Proposition}
Let $I$ be a homogeneous ideal in $R$ and $\Omega$ a truncating
subset. Let $C_\Omega[\omega_1]$ and 
$C_\Omega[\omega_2]$ be maximal cells in the truncated 
Gr\"obner fan $\cF_\Omega(I)$ of $I$. Suppose that $G$ is the
reduced $\Omega$-Gr\"obner basis for $I$ over $\omega_1$. 
If $\omega\in \overline{C_\Omega[\omega_1]} \cap 
\overline{C_\Omega[\omega_2]}$, then
\begin{enumerate}[(i)]
\item
The reduced $\Omega$-Gr\"obner basis for $\ini_\omega(I)$
over $\omega_1$ is $G_\omega = \{\ini_\omega(g) \mid
g\in G\}$.
\item
If $H$ is the reduced $\Omega$-Gr\"obner basis for $\ini_\omega(I)$
over $\omega_2$, then
$$
\{f - f^G \mid f\in H\}
$$
is a minimal $\Omega$-Gr\"obner basis for $I$ over $\omega_2$ 
modified by $\omega$.
\item
The reduced $\Omega$-Gr\"obner basis for $I$ over $\omega_2$ modified by
$\omega$ coincides with the reduced $\Omega$-Gr\"obner basis for $I$ over
$\omega_2$. 

\end{enumerate}
\end{Proposition}
\begin{Proof}
The items (i) and (iii) follow as in  Proposition 3.2 of \cite{FJLT} taking
Proposition \ref{PropositionRG} and Corollary \ref{CorollaryDeformed} into
account. For the proof of (iii) observe that if $H' = \{f_1, \dots, 
f_s\}$ is the reduced Gr\"obner basis of $\ini_\omega(I)$ over
$\omega_2$, then $G'' := \{f_1 - f_1^{G'}, \dots, f_s - f_s^{G'}\}$ is a
minimal Gr\"obner basis for $I$ over $\omega_2$ modified by $\omega$ by
Proposition 3.2(ii) in \cite{FJLT}. Here $G'$ is the reduced Gr\"obner
basis for $I$ over $\omega_1$. Notice that $G''$ consists of homogeneous
elements and that 
$$
G''_\Omega = \{g\in G'' \mid \deg(g) \in \Omega\}
$$ 
is a minimal $\Omega$-Gr\"obner basis for $I$ over $\omega_2$ modified
by $\omega$. If $\deg(f_i - f_i^{G'})\in \Omega$, then 
$\deg(f_i)\in \Omega$. In this case $f_i^{G'} = f_i^G$. Using
Proposition \ref{PropositionRG} this
finishes the proof of (ii).
\end{Proof}

Now the generic Gr\"obner
walk (\cite{FJLT}, \S 4) carries over verbatim to the truncated setting using
the Buchberger algorithm with truncation in step (iv).

\section{Lattice ideals}

In the rest of this paper we will remain exclusively in the setting
of lattice ideals. 
Recall the decomposition of an integral vector 
$v\in \ZZ^n$ into $v = v^+ - v^-$,
where $v^+, v^-\in \NN^n$ are vectors with disjoint support.
For $u, v\in \NN^n$ we let $u \leq v$ denote the partial
order given by $v - u\in \NN^n$.
For a subset
$\cB \subset \ZZ^n$ we associate the ideal
$$
I_\cB = \<x^{v^+} - x^{v^-} \mid v\in \cB\> \subset R.
$$
In the case where $\cB = \cL$ is a lattice we call
$I_\cL$ the {\it lattice ideal\/} associated to $\cL$. 
Recall that lattice ideals are saturated i.e. if $f\in I_\cL$
is divisible by a variable $x_i$, then $f/x_i \in I_\cL$.
This means that we apply the homogeneous Buchberger algorithm
with sat-reduction as explained in \cite{Lau}.

Define 
$$
\bin(w) = x^{w^+} - x^{w^-}
$$
for $w\in \ZZ^n$.
The (saturated) $S$-polynomial of $\bin(u)$ and $\bin(v)$ is
then given by $\bin(u-v)$. Similarly if 
$v^+ \leq u^+$ we may reduce $\bin(u)$ by
$\bin(v)$ giving $\bin(u-v)$. We have silently assumed that
the initial term of $\bin(w)$ is $x^{w^+}$ for the term ordering
in question.

Usually a generating set $\cB$ for $\cL$ as an abelian group is given. 
Computing the lattice ideal $I_\cL \supset I_\cB$ can be done using that
$$
I_\cL = I_\cB : (x_1 \cdots x_n)^\infty,
$$
where
$I:f^\infty$ denotes the ideal given by
$$
\{r\in R \mid r\, f^m \in I, \mbox{\ for\ } m \gg 0\}
$$
for an ideal $I\subset R$ and an element $f\in R$ (\cite{St}, Lemma 12.2).
If $\cB$ contains a positive vector, then $I_\cB = I_\cL$ (\cite{St}, 
Lemma 12.4). If $\cL \cap \NN^n = \{0\}$, 
$I_\cL$ may be computed from $I_\cB$ using 
Gr\"obner basis computations for different reverse 
lexicographic term orderings (\cite{St}, Lemma 12.1).

\subsection{The generic Gr\"obner walk for lattice ideals}

We now specialize the generic Gr\"obner walk to the setting
of lattice ideals representing binomials by integer
vectors as $\bin(v)$.
In our implementation of the generic walk we walk from
$\prec_{c_1}$ to $\prec_{c_2}$, where
$\prec$ is the reverse lexicographic order given by
$x_1 \prec \cdots \prec x_n$ and $c_1, c_2$ are integer
vectors.

In the algorithm outlined below we walk between two arbitrary
multiplicative orderings $\prec_1, \prec_2$ and move
inside Gr\"obner cones given by minimal Gr\"obner bases (cf.~\cite{FJLT}, 
Proposition 2.3). Autoreduction is replaced by a simplified
lifting step.

In the notation of (\cite{FJLT}, \S2.3) we have
$$
\delta_\prec(\bin(w)) = w
$$
assuming that $x^{w^-} \prec x^{w^+}$.
The {\it facet preorder\/} $\prec$ is now given on binomials 
$\bin(u)$ and $\bin(v)$ as in (\cite{FJLT}, \S4, (3)) by
$$
\bin(u) \prec \bin(v) \iff T u v^t \prec_1 T v u^t,
$$  
where $T$ is a matrix defining the target term order $\prec_2$. 
We get as in (\cite{FJLT}, \S4) that $\bin(u) \prec 
\bin(v)$ and $\bin(v) \prec \bin(u)$ imply that $u$ is a 
multiple of $v$. If $G := \{\bin(v_1), \dots, \bin(v_r)\}$ is a
Gr\"obner basis then $v_1, \dots, v_r$ lie in a common half space and
the facet preorder
induces a total ordering on $G$. In particular one gets that the
facet in the generic Gr\"obner walk is given by a unique {\it
facet binomial\/}. This means that Gr\"obner basis computations
on facets proceed as in (\cite{HT}, Algorithm 3.1). To give some
more details we introduce the notation
$$
\mon(w) = x^{w^+}
$$
for $w\in \ZZ^n$.
Suppose that $G$ above is a minimal Gr\"obner basis and that
$\bin(v_1)$ is minimal in the facet preorder. Then step (c)
of (\cite{HT}, Algorithm 3.1) is to compute a Gr\"obner basis
of the ideal 
\begin{equation}\label{facet}
\<\bin(-v_1), \mon(v_2), \dots, \mon(v_r)\>.
\end{equation}
Working with minimal Gr\"obner bases it may happen that $\mon(-v_1)$
is divisible by a  monomial $\mon(v_j), j = 2, \dots, r$. In this
case the reduced Gr\"obner basis of (\ref{facet}) is
$$
\{\mon(v_1), \dots, \mon(v_r)\}
$$
and $\mon(v_1)$ lifts to $\bin(w)$, where $\bin(-w)$ is the reduction
of $\bin(-v_1)$ modulo $$
\{\bin(v_2), \dots, \bin(v_r)\}.
$$ 
In this
way the usual autoreduction of the Gr\"obner walk is built into
the lifting. Notice that $\mon(v_j)$ lifts to $\bin(v_j)$ for
$j = 2, \dots, r$. These observations account for step (ii) 
in {\bf facet\_buchberger} below. On the other hand, if 
$\mon(-v_1)$ is not divisible by any of the monomials, then
$\bin(v_1)$ is a (real) facet binomial of a facet in the Gr\"obner
cone corresponding to the reduced Gr\"obner basis.
In this case we end up with a minimal Gr\"obner basis
$$
\{\bin(-v_1), \mon(w_1), \dots, \mon(w_s)\}
$$
of the ideal in (\ref{facet}). This lifts 
to the minimal Gr\"obner basis
$$
\{\bin(-v_1), \bin(w_1), \dots, \bin(w_s)\}.
$$
The details of the algorithm are given below. 
The variable {\it facet\_list\/} 
contains a list of binomials ordered in ascending order according to
the facet preorder (these are potential facet binomials). The variable $G$
contains the current minimal Gr\"obner basis. The procedure
{\bf initialize\_facet\_list} initializes $G$ and {\it facet\_list\/} given
$B$. The procedure {\bf insert} inserts a given binomial into $G$ and
updates {\it facet\_list\/}. 

Notice that we do not really compute the $S$-polynomials in (iv.b) below.
We optimize the algorithm by replacing the $S$-polynomial
$S(\bin(w), \mon(v))$ with the initial term $x^{(w - v)^+}$ of
the saturated $S$-polynomial $\sat(S(\bin(w), \bin(v))$.

\begin{Algorithm}[Generic Gr\"obner walk for lattice ideals]\label{AlgorithmGW}
\

\

\noindent
{\bf INPUT}: Integer vectors $c_1, c_2$. Integer vectors 
$B = \{v_1, \dots, v_r\}$ such that 
$\{\bin(v_1), \dots, \bin(v_r)\}$ is a minimal Gr\"obner basis for
$I_\cL$ over $\prec_1$.

\

\noindent
{\bf OUTPUT}: Integer vectors $G = \{w_1, \dots, w_s\}$ such
that $\{\bin(w_1), \dots, \bin(w_s)\}$ is a minimal Gr\"obner
basis over $\prec_2$.

\begin{enumerate}[(i)]
\item {\bf initialize\_facet\_list};
\item while $(facet\_list \neq \emptyset)$ do
\begin{enumerate}
\item $facet\_bin := \mbox{first element in $facet\_list$}$
\item {\bf facet\_buchberger};
\end{enumerate}
\end{enumerate}

\

\noindent
{\bf facet\_buchberger:}
\begin{enumerate}[(i)]
\item
Delete $facet\_bin$ from $G$ and put $bin := -facet\_bin$;
\item
if $w^+\leq bin^+$ for some $w\in G$

{\narrower \parindent 25 true pt
  reduce $bin$ by $G$;

{\bf insert}(-bin);

return;
}
\item
Spairs := $\emptyset$;
\item
for $v$ in $G$ do
\begin{enumerate}
\item
if $(bin^+ \wedge v^+ = 0$)

{\narrower \parindent 25 true pt 
  continue;
}

\item
$Spairs := Spairs \cup \{v - bin\}$
\end{enumerate}

\item
Delete $v\in G$ if $bin^+\leq v^+$;

\item
while ($Spairs \neq \emptyset$) do

\begin{enumerate}
\item
Select $s$ in $Spairs$ and put $Spairs := Spairs \setminus \{s\}$.
\item
Reduce $s$ by $bin$;
\item
if $(v^+ \leq s^+)$ for some $v\in G$

{\narrower \parindent 25 true pt
  continue;
}

\item
Delete $v\in G$ if $s^+ \leq v^+$.

\item
{\bf insert}(s)



\item
$Spairs := Spairs\cup \{s - bin\}$

\end{enumerate}

\item
{\bf insert}(bin);

\end{enumerate}

\end{Algorithm}

{\it Truncation\/} blends in easily with Algorithm \ref{AlgorithmGW}. Suppose 
that $\Omega$ denotes the truncating subset. First binomials in $B$ with
degrees outside $\Omega$ are discarded.
With every addition of an $S$-binomial in step iv(b) of 
{\bf facet\_buchberger}, a test for degree membership of the truncating
subset $\Omega$ is done. If the test
fails for the $S$-binomial it is not added to $Spairs$.

\begin{Example}
We give a very simple example illustrating Algorithm \ref{AlgorithmGW}. Consider the
ideal
$$
I = \<x - t^2, y - t^3\> \subset k[t, x, y].
$$
Clearly $G := \{x - t^2, y - t^3\}$ is a Gr\"obner basis for $I$ over
the weight vector $\sigma = (-1, 0, 0)$, where $t$ corresponds to $(1, 0, 0)$ etc.
We wish to walk to the weight vector $\tau = (1, 0, 0)$ breaking ties with
the reverse lexicographic order given by $t <  x < y$. Let $\prec$ denote
the corresponding facet preorder. Then $x - t^2 \prec y - t^3$ and we
begin by ``computing'' a Gr\"obner basis for $\<t^2 - x, y\>$ giving
$G = \{t^2 - x, y - t^3\}$ after lifting. In the
following step the facet binomial is $y - t^3$, which gets replaced by
the reduction $y - t x$ in step (ii) of {\bf facet\_buchberger}. This
accounts for the next facet binomial. We then compute a Gr\"obner basis
of $\<t x - y, t^2\>$ giving $\{t x - y, t^2, t y, y^2\}$. This lifts
to $\{t x - y, t^2 - x, t y - x^2, y^2 - x^3\}$, which is the reduced Gr\"obner
basis for $I$ over $(1, 0, 0)$, since $t y - x^2$ and $y^2 - x^3$ are
not candidates for facet binomials as the vectors $(1, -2, 1)$ and
$(0,  -3, 2)$ both are outside $C_{<_\sigma, <_\tau}$ (cf.~\S4 of
\cite{FJLT}).
\end{Example}

\section{Computational experience}

In \cite{AL} a collection of integer knapsacks are constructed
related to the classical Frobenius problem of finding the
largest number, which is not a sum of given relatively prime
natural numbers. Feasibility for these knapsacks turn out to 
be very hard for traditional branch and bound software like CPLEX, but
easy for lattice reduction methods as shown in \cite{AL}.

In \cite{LATTE} these knapsacks are equipped with a feasible
right hand side and a specific cost vector $c$. 
These examples form
the point of departure in this section, where we specifically
document performance for computing (truncated) test sets
using the package {\tt GLATWALK}\footnote{\tt home.imf.au.dk/niels/GLATWALK}. It turns out that test sets
in the feasibility case is by far the hardest computations. Test
sets with respect to the cost vector in \cite{LATTE} finish
in negligible timings ($< 0.05$ seconds) using both the generic
walk and the Buchberger algorithm with truncation.

Each of the examples are of the form:
maximize $c x$, where
\begin{equation}\tag{$*$}
A x = b,
\end{equation}
$x\in \NN^n$ and $A$ is a $1\times n$-matrix $(a_1 \cdots a_n)$. 
The cost vector $c$ and the matrices $A$ may be found in \cite{LATTE}.
The first 
step is finding a feasible 
solution to ($*$). As in \cite{FJLT} this results in the
knapsack: minimize $t$ subject to $A x + t = b$, where
$t\in \NN$ and $x\in \NN^n$. This leads to the problem
of finding a (truncated) Gr\"obner basis of
\begin{equation}\tag{$**$}
\<x_1 - t^{a_1}, \dots, x_n - t^{a_n}\>
\end{equation}
with respect to the vector $\tau = (1, 0, \dots, 0)$, where $t$
is the ``first'' variable. We may compute
this Gr\"obner basis directly using the Buchberger algorithm or walk
from the vector $\sigma = (-1, 0, \dots, 0)$.
The performance of the functions 
{\bf walk}  and {\bf gbasis} of
{\tt GLATWALK} 
for computing a full Gr\"obner 
basis of ($**$) over $\tau$
is reported in \cite{FJLT}.
The second step is the computation of the toric ideal $I_A$ (associated
with the integer matrix $A$) and
its Gr\"obner basis over the vector $-c$. The function
{\bf saturate} of {\tt GLATWALK} performs the saturation
necessary in computing lattice ideals. Below\footnote{All timings
are in seconds. The computations were carried out
on an ACER notebook 1.6 GHz Pentium mobile with 1MB L2 cache.}
 we have
computed the ideals $I_A$ using {\bf saturate} after 
LLL-reducing $\mathrm{Ker}(A)$ with the function {\bf LLL}. 
In many of the examples, LLL-reduction offers great savings in 
the computation of the saturation.
The third column shows the timing of {\bf
gbasis} in computing a full Gr\"obner basis over
$-c$ for $I_A$. The fourth column is the timing of {\bf walk} in walking
from $-e_1$ to $-c$. The 
fifth and sixth columns show sizes of
the full and truncated reduced Gr\"obner bases of $I_A$ over $-c$.

\

\

\

$$
\begin{tabular}{|p{1.8cm}|r|r|r|r|r|r|}
\hline
EXAMPLE & {\bf saturate} & {\bf gbasis} & {\bf walk} &size & tr size \\
\hline
cuww1 & 0.1 &   3.1 & 1.8 & 2618 & 7\\
\hline
cuww2 & 0.0 &  0.3 & 1.3 & 898 &16\\
\hline
cuww3 & 0.1 &  0.4 & 10.1 & 963 & 16\\
\hline
cuww4 & 0.0  & 3.3 & 59.6 & 3143 & 5\\
\hline
cuww5 & 0.0  & 0.0 & 102.8 & 267 & 32\\
\hline
prob1 & 0.0  & 0.0 & 2.7 & 180 & 75\\
\hline
prob2 & 0.0  & 0.0 & 0.5 & 280 & 45\\
\hline
prob3 & 0.5  & 0.0 & 2.5 & 163 &94\\
\hline
prob4 & 0.2  & 0.1 & 122.3 & 475 & 83\\
\hline
prob5 & 0.6  & 0.0 & 0.0& 68 & 56\\
\hline
prob6 & 0.4  & 9.2 & 39.5 & 4541 & 94\\
\hline 
prob7 & 0.4  & 1.8 & 72.8 & 2036 & 79\\
\hline
prob8 & 0.9 & 0.0 & 2.4 & 227 & 103\\
\hline
prob9 & 0.0 &  0.0 & 0.5 & 108 & 108\\
\hline
prob10 & 1.4 & 0.1 &  517.5 & 536 & 119\\
\hline
\end{tabular}
$$

\ 

\

\

Both {\bf walk} and {\bf gbasis} finish in negligible timings ($<0.05$
seconds) in computing the {\it truncated\/} Gr\"obner bases. However in computing
the full Gr\"obner bases in the above table, {\bf walk}  does not compare well 
with {\bf gbasis}. Typically to compute a target 
Gr\"obner basis with less than $1,000$ binomials, the generic walk 
traverses cones associated with reduced Gr\"obner bases of more 
than $20,000$ binomials along a straight line intersecting many
cones in the Gr\"obner fan.

Most of the examples above indicate that the truncated Gr\"obner fan
is much smaller than the full Gr\"obner fan. The straight line path
in the truncated Gr\"obner fan traverses significantly fewer cones. It is open
for further research exactly when the walk is a substantial improvement (as
in many of the feasibility examples reported in \cite{FJLT}). Perhaps
a combination of direct Gr\"obner basis computations for suitably 
chosen (easier) weight vectors tending to the target vector followed by a 
walk to the target order may lead to improvements.

\end{document}